\documentclass[11pt]{amsart}
\usepackage{amssymb,amsmath, amsthm, amsfonts, float, mathtools, ytableau, xcolor}
\usepackage[colorlinks=true, urlcolor=black,linkcolor=blue, citecolor=blue]{hyperref}
\usepackage{bookmark}
\DeclarePairedDelimiter{\ceil}{\lceil}{\rceil}
\DeclarePairedDelimiter{\floor}{\lfloor}{\rfloor}

\newtheorem{theorem}{Theorem}[section]

\newtheorem{proposition}{Proposition}[section]
\newtheorem{lemma}{Lemma}[section]
\newtheorem{definition}{Definition}[section]

\newtheorem{remark}{Remark}[section]
\newcommand{\bbox}{\hfill $\Box$}
\newcommand{\pf}{\noindent {\it Proof:} }

\begin{document}

\title[\resizebox{4.5in}{!}{The random ({\large\MakeLowercase{n-k}})-cycle to transpositions walk on the symmetric group}]{The random {\Large\MakeLowercase{(\textit{n-k})}}-cycle to transpositions walk on the symmetric group} 

\author{%
Alperen Y. \"{O}zdemir}

\address{%
University of Southern California\\
Department of Mathematics\\
Los Angeles, California, 90089-2532\\
Tel: +1 (323) 596-5095 \\
E-mail: aozdemir@usc.edu \\
ORCID ID: 0000-0003-2730-7240}



\begin{abstract}
We study the rate of convergence of the Markov chain on $S_n$ which starts with a random $(n-k)$-cycle for a fixed $k \geq 1$, followed by random transpositions. The  convergence to the stationary distribution turns out to be of order $n$. We show that after $cn + \frac{\ln k}{2}n$ steps for $c>0$, the law of the Markov chain is close to the uniform distribution. The character of the defining representation is used as test function to obtain a lower bound for the total variation distance. We identify the asymptotic distribution of the test function given the law of the Markov chain for the $(n-1)$-cycle case. The upper bound relies on estimates for the difference of normalized characters. 
\end{abstract}

\maketitle

\noindent
\textbf{Acknowledgement}
The author would like to thank Jason Fulman for suggesting the problem and his most valuable comments.

\section{Introduction} \label{intro}

Random walks on symmetric group are widely associated with card shuffling problem. For an extensive survey on random walks defined on finite groups, see \cite{D} or \cite{S}. In their influential work \cite{DS}, Diaconis and Shahshahani study the random transposition walk on symmetric group by using representation theory. They demonstrate that the chain exhibits a cutoff at $\frac{1}{2}n \ln n$ for the total variation distance. In recent years, more probabilistic approaches are developed. Berestycki et al. \cite{BSZ} prove a cutoff at $\frac{n}{k} \ln n$ steps for $k-$cycle walk by analyzing the evolution of the cycle distribution. A path coupling argument is shown to be applicable for random transposition walk to obtain a mixing time of order $n \ln n$ in \cite{Bo}.  \\

We study a variation of the random transposition walk, following the techniques developed in \cite{DS}, by taking the initial permutation to be a cycle with a fixed number of fixed points instead of a transposition. The upper bound for the $(n-1)-$cycle to transpositions case is given in \cite{DP} along with both an upper and a lower bound for $n$-cycle to transpositions walk. The lower bound is left an open problem.  \\

The chain starts at the identity. At the initial step, an $(n-k)$-cycle is uniformly selected for a fixed $k \geq 1.$ From then on, the transition is selecting a transposition uniformly and multiplying it by the permutation of the current state to obtain another permutation for the next state. Observe that the chain alternates between $A_n$ and $\medmuskip=0mu S_n \setminus A_n$. Therefore the limiting distribution is not the uniform distribution on $S_n$ unlike the random transposition walk in \cite{DS}, for which they define a holding probability that makes the random walk lazy. We define the transition probabilities as in \cite{DP}, and allow the limiting distribution to be different for even steps and odd steps. The analysis is technically the same and the limiting distribution is uniform on the corresponding coset. $U_t$ will stand for the uniform distribution on the coset ($A_n$ or $ \medmuskip=0mu S_n \setminus A_n$) where the chain is on, at step $t$. \\

To be more precise on the stationary distribution, we identify it with respect to the length of the initial cycle, $n-k.$ We call two integers to have the same parity if they are both odd or both even, otherwise we call them to have different parity. Denoting the uniform distribution over a set $S$ by $\mathcal{U}(S),$ we have 
\[U_t= \begin{cases} 
      \mathcal{U}(A_n) & \text{if $t$ and $n-k$ have the same parity,} \\
      \mathcal{U}(\medmuskip=0mu S_n \setminus A_n) & \text{if $t$ and $n-k$ have different parity.} 
   \end{cases}
\]

Let the probability assigned to  $\sigma \in S_n$  be $\mu_{t}(\sigma)$ after $t$ steps. The total variation distance between $\mu_{t}$ and the uniform distribution $U_{t}$ is

\begin{equation} \label{TV}
\| \mu_{t} - U_{t} \|_{TV} = \frac{1}{2} \sum_{\sigma \in S_n} |\mu_{t}(\sigma) - U_{t}(\sigma)| = \max_{S \subseteq S_n}| \mu_{t}(S) - U_{t}(S)|.
\end{equation}
We are interested in finding upper and lower bounds for the probability metric defined above and identify the convergence rate. For the asymptotic convergence, we use the notation  $f(n) \sim g(n),$ meaning $\lim_{n \rightarrow \infty} \frac{f(n)}{g(n)}=1,$ where $n$ is the number of elements that $S_n$ is defined over.

The organization of the paper is as follows. In Section \ref{chp1} we summarize the group representation theory techniques that are used in the rest of the paper. In Section \ref{chp2}, we focus on the $(n-1)$-cycle case. The moments of the character of the defining representation are found and are shown to determine a unique distribution. We evaluate the distribution function to bound the total variation distance from below. Combining with Ding's result \cite{DP} we obtain 
\begin{theorem} \label{mthm1}
As $n \rightarrow \infty$, for any $c > 0$, after a random $(n-1)$-cycle and $cn$ random transpositions,
\begin{equation*}
\frac{ 1-e^{e^{-2c}} + e^{e^{-2c}-2c} }{e} - o(1) \leq \| \mu_{cn+1} - U_{cn+1} \|_{TV} \leq \frac{e^{-4c}}{2 \sqrt{1-e^{-4c}}} + o(1).
 \end{equation*} 
\end{theorem}

In Section \ref{chp3}, a lower bound on the total variation distance for $k \geq 2$ after $cn + \frac{\ln k}{2}n$ steps is obtained. The first two moments  of the character of the defining representation are needed to bound the distance to the stationary distribution. Then we derive estimates for the differences of characters and combine them with the previous results in literature to find the upper bound for the given convergence time. This leads to the following theorem.
\begin{theorem} \label{mthm2}
Let $k\geq 2$ be a fixed number. As $n \rightarrow \infty$, for any $c > 0,$ after a random $(n-k)$-cycle and random $t=cn + \frac{n}{2} \ln{k}$ transpositions, 
\begin{equation*}
\frac{1}{12}\, e^{-4c} - o(1) \leq \| \mu_{t+1} - U_{t+1} \|_{TV} \leq \sqrt{\frac{e-1}{2}} \, e^{-2c} +o(1).
\end{equation*}
\end{theorem}

\section{Representation theory techniques} \label{chp1}

We summarize the techniques that are employed in the following sections. The tools developed in \cite{DS} to study the random transposition walk are used to obtain lower and upper bounds in our case. To start with the connection to the group representation theory, first consider the \textit{Fourier transform} of the measure $\mu$ defined on a finite group $G$, evaluated at any representation $\lambda$ of $G$, 

\begin{equation}\label{fourier}
\widehat{\mu}(\lambda)=\sum\limits_{g \in G} \mu(g) \lambda(g) 
\end{equation}
The inverse transform is expressed as
\begin{equation*}
\mu(g)= \sum\limits_{ \lambda \, irred.}d_{\lambda} tr(\lambda(g^{-1})\widehat{\mu}(\lambda))
\end{equation*}
where the sum is over all irreducible representations of $G$ and $d_{\lambda}$ stands for the \textit{dimension of the representation} $\lambda,$ which is the number of standard Young tableaux of shape $\lambda.$ The inverse transform leads to 
Plancherel's formula stated below.
\begin{equation} \label{planc}
\sum\limits_{g \in G} |\mu(g)|^2 = \frac{1}{|G|} \sum\limits_{ \lambda \, irred.}d_{\lambda} tr(\widehat{\mu}(\lambda)\widehat{\mu}(\lambda)^{T}).
\end{equation}
This establishes the connection with the total variation distance if we take $\mu$ on the left hand side to be the difference of two measures.\\

On the right hand side, we have the trace of Fourier transforms, which can be evaluated by writing it in terms of the characters of the irreducible representations. We denote the character of the representation $\lambda$ by $\chi_{\lambda},$ which is defined as $\chi_{\lambda}(g)=tr(\lambda(g)).$ An important fact about characters is that they are class measures, so they are constant on conjugacy classes. We have the following lemma, which can be found in Chapter 16 of \cite{B}.

\begin{lemma} \emph{(\cite{B})} \label{trace} Let a class measure $\mu$ be defined on $G$ and $\lambda$ be a representation of $G.$ Then,
\begin{equation*}
\widehat{\mu}(\lambda)= \frac{1}{d_{\lambda}} \left[\sum\limits_{g \in G} \mu(g) \chi_{\lambda}(g) \right] I_{d_{\lambda}}. 
\end{equation*}
\end{lemma}

We also note a fact about the Fourier transform of the convolution of measures, as a Markov chain at step $t$ can be viewed as  $t-$fold convolution measure of the transition probabilities. The convolution of two measures $\mu$ and $\nu$ is defined to be
\begin{equation*}
\mu * \nu (g) = \sum\limits_{h \in G} \mu(g h^{-1}) \nu(h) 
\end{equation*}
The Fourier transform of a convolution satisfies 
\begin{equation} \label{conv}
\widehat{\mu * \nu}(g)= \widehat{\mu}(g) \widehat{\nu}(g). 
\end{equation}

\subsection{Lower bound}

It follows from the definition of the total variation distance (\ref{TV}), after a random $(n-1)$-cycle and $cn$ transpositions, we have
\begin{equation} \label{lb}
\| \mu_{cn+1} - U_{cn+1} \|_{TV} \geq |\mu_{cn+1}(A) - U_{cn+1}(A)|
\end{equation}
where $A$ is the set of fixed point free permutations. This particular choice of the subset of $S_n$ establishes a connection with the symmetric group representations.

We first introduce the \textit{defining representation} of $S_n$, which is the $n$-dimensional representation $\rho$ where 
\begin{equation*}
\rho(\sigma)(i,j)=
 \left\{
\begin{array}{ll}
      1& \sigma(j)=i \\
      0 & \textrm{otherwise} \\
\end{array} 
\right. 
\end{equation*}
for $\sigma \in S_n.$ \\

Denote the character of the representation, or the traces of the matrices given above, by $\chi_{\rho}.$ So that $\chi_{\rho}(\sigma)$ counts the number of fixed points of $\sigma \in S_n.$ Therefore $\sigma \in A$ if and only if $\chi_{\rho}(\sigma)=0.$\\

The number of fixed points of a uniformly random permutation has limiting Poisson distribution with parameter $1$; elementary proofs of this fact can be found in \cite{T}. The limiting distribution in our case, which alternates between the uniform distribution over permutations in $A_n$ and in $\medmuskip=0mu S_n \setminus A_n$, is also asymptotically $\mathcal{P}(1)$. A proof is given in \cite{DP} by showing that the $r^{th}$ moment of $\chi_{\rho}$ under $U_t$ agrees with the $r^{th}$ moment of $\mathcal{P}(1)$ for $r \leq n,$ i,e.,  
\begin{equation*}
\mathbf{E}_{U_t}(\chi_{\rho}^r)=\sum\limits_{i=0}^r {r \brace i}=B_r
\end{equation*}
for $r \leq n,$ where ${r \brace i}$ is a notation for the Stirling numbers of second kind. The sum above gives the $r^{th}$ Bell number. \\

Next we consider the moments of $\chi_{\rho}$ under the law of the Markov chain at step $t.$ The defining representation is reducible and decomposed as

\begin{equation} \label{fmom}
\rho=S^{(n)} \oplus S^{(n-1,1)}
\end{equation}
where $S^{\lambda}$ is the Specht module associated with partition $\lambda$ of $n$, noting that the partitions of $n$ are in one to one correspondance with the irreducible representations of $S_n.$ See Chapter 4 of \cite{J} for details.

In order to find the higher moments, we consider the decomposition of the tensor product $\rho^{\otimes r}$.
\begin{equation} \label{rmom}
\rho^{\otimes r}=\underset{\lambda \vdash n}{\oplus} a_{\lambda,r} S^{\lambda}
\end{equation}
where $\lambda \vdash n$ means that $\lambda$ is a partition of $n$. Then we use the facts below,
\begin{equation*}
 \begin{split}
\chi_{\rho_1 \oplus \rho_2}&=\chi_{\rho_1}+ \chi_{\rho_2}, \\
\chi_{\rho_1 \otimes \rho_2}&=\chi_{\rho_1} \cdot \chi_{\rho_2},
 \end{split}
\end{equation*}
and invoke some Fourier analytic results, which can be found in \cite{B} Chapter 16 with detailed proofs and in \cite{D} Section 2C in most relevance to our case, to calculate the moments needed.
\begin{equation}\label{di}
\mathbf{E}_{\mu}((\chi_{\rho})^r)=\sum_{\lambda \vdash n } a_{\lambda,r} tr(\hat{\mu}(\lambda)).
\end{equation}

Now define $\bar{\lambda} = (\lambda_2,\lambda_3,...,\lambda_m)$, the partition obtained by removing the first row of $\lambda$. The generating function of coefficients ${a_{\lambda,r}}$ for a fixed partition $\lambda$, 
\begin{equation*}
\sum\limits_{r \geq |\bar{\lambda}|} {a_{\lambda,r}} \frac{x^r}{r!}= \frac{d_{\lambda}}{|\bar{\lambda}|!}e^{e^{x}-1}(e^x-1)^|\bar{\lambda}|,
\end{equation*}
is found by Goupil and Chauve in \cite{GC}. The coefficients are identified by Ding in \cite{DP} to be
\begin{equation} \label{coeff}
a_{\lambda,r}=d_{\bar{\lambda}}\sum\limits_{i=|\bar{\lambda}|}^r \binom{i}{|\hat{\lambda}|}{r \brace i}.
\end{equation}
for $1 \leq r \leq n- \lambda_2.$ Therefore, explicit expressions for higher moments of $\chi_{\rho}$ can be obtained by finding the traces of the Fourier transforms of $\mu.$ Lemma \ref{trace} and the fact about the convolution of measures (\ref{conv}) are applied in that respect. Particularly in our case,
 \begin{equation*}
\text{tr}(\hat{\mu_t}(\lambda))= d_{\lambda}\left( \frac{\chi^{\lambda}_{(n-k,1^{k})}}{d_{\lambda}} \right)\left( \frac{\chi^{\lambda}_{(2,1^{n-2})}}{d_{\lambda}} \right)^t 
\end{equation*}
where $\chi^{\lambda}_{(n-k,1^k)}$ is the character of the representation $\lambda$ evaluated at an $(n-k)-$cycle. \\

Very often it might be difficult to calculate the the higher moments for an arbitrary permutation. The following lower bound lemma requires the first moment and a higher one to provide a lower bound on the total variation distance. We follow the proof of a result in Chapter 7 of \cite{L}, which gives a lower bound using the first two moments. 
\begin{lemma} \label{p2}
Let $\mu$ and $\nu$ be two probability distributions on $\Omega$, and $X$ be a random variable from $\Omega$ to a finite subset $A$ of $\mathbb{R}$. Then for $p,q \in (1, \infty)$ satisfying $\frac{1}{p}+ \frac{1}{q}=1$,
\begin{equation*} 
  \| \mu - \nu \|_{TV} \geq \frac{1}{2} \, \frac{[\mathbf{E}_{\mu}(X)-\mathbf{E}_{\nu}(X)]^q}{\left(\mathbf{E}_{\mu}(X^p)+\mathbf{E}_{\nu}(X^p)\right)^{q-1}}.
\end{equation*}

\end{lemma}
\pf
We first define the average measure on $A$,
 \begin{equation*}
 \sigma(x)= \frac{P_{\mu}(X=x) + P_{\nu}(X=x)}{2},
 \end{equation*}
  and 
the functions
\begin{equation*}
\alpha(x)=\frac{P_{\mu}(X=x)}{\sigma(x)}, \quad \quad \beta(x)=\frac{P_{\nu}(X=x)}{\sigma(x)}.
\end{equation*}
Then by H\"{o}lder's inequality,
\begin{equation} \label{hol}
\begin{split}
\mathbf{E}_{\mu}(X)-\mathbf{E}_{\nu}(X)&=\sum\limits_{x\in A} x(\alpha(x)-\beta(x)) \sigma(x)\\
 &\leq \left( \sum\limits_{x\in A} x^p\sigma(x) \right)^{1/p} \left(\sum\limits_{x \in A} |\alpha(x)-\beta(x)|^q \sigma(x)\right)^{1/q} 
\end{split}
\end{equation}
where $\frac{1}{p}+ \frac{1}{q}=1$ for $p,q \in (1, \infty).$ Observe that
\begin{equation*}
\sum\limits_{x \in A} x^p \sigma(x) = \frac{\mathbf{E}_{\mu}(X^p)+\mathbf{E}_{\nu}(X^p)}{2} , \quad \quad \sum\limits_{x \in A} |\alpha(x)-\beta(x)| \sigma(x)= 2\| \mu - \nu \|_{TV}
\end{equation*}
and
\begin{equation*}
|\alpha(x)-\beta(x)|=\frac{2|P_{\mu}(X=x)-P_{\nu}(X=x)|}{P_{\mu}(X=x)+P_{\nu}(X=x)} \leq 2.
\end{equation*}
Combining the observations above in $\eqref{hol}$, we arrive at the result.

\bbox

\subsection{Upper bound}
An upper bound lemma is provided in Chapter 3B of \cite{D}. Yet we need a slight modification of it due to the alternation of the stationary distribution $U_t$ between $A_n$ and $\medmuskip=0mu S_n \setminus A_n$ in our case. \\

We first introduce the notation for the two one dimensional irreducible representations of $S_n.$ The first one is $\lambda_{triv}$ is the trivial representation that maps every element to $1$. The other one is $\lambda_{sign}$, which is $1$ for even permutations and $-1$ for odd permutations. 
The upper bound lemma is as follows.
\begin{lemma}
Let $\mu_t$ be the law of the Markov chain at step $t$ and $U_t$ be the stationary distribution described in (\ref{intro}). Then,
\begin{equation*}
4\| \mu_{t} - U_{t} \|_{TV}^2 \leq \, \frac{1}{2} \sum_{\substack{\lambda \, irred. \\ \lambda \neq \lambda_{triv}, \lambda_{sign}}} d_{\lambda} tr(\widehat{\mu_t}(\lambda)\widehat{\mu_t}(\lambda)^T).
\end{equation*}
\end{lemma}
\pf We first argue for the following facts on the Fourier transform of $U_t$. 
\begin{equation}\label{triv}
\begin{split}
\widehat{U_t}(\lambda) &=\widehat{\mu_t}(\lambda) \, \textrm{ if } \lambda= \lambda_{triv} \text{ or } \lambda_{sign}, \\
\widehat{U_t}(\lambda) &\text{ is the zero matrix otherwise.} \\
\end{split}
\end{equation} 
Clearly, $\widehat{U_t}(\lambda_{triv})=1.$ $\widehat{U_t}(\lambda_{sign})$ is either $1$ or $-1$ depending on the coset the chain is on at step $t.$ Suppose that the chain is on $A_n.$ Then
\begin{equation*}
\widehat{\mu_t}(\lambda_{sign})= \sum_{\sigma \in A_n} \mu_t(\sigma) \chi_{\lambda_{sign}}(\sigma) = \sum_{\sigma \in  A_n} \mu_t(\sigma) = 1,
\end{equation*}
which is equal to $\widehat{U_t}(\lambda_{sign})$. If the chain on $\medmuskip=0mu S_n \setminus A_n,$ they are both $-1.$ If $\lambda$ is different from the trivial and the sign representation, $\widehat{U_t}(\lambda)$ is the zero transformation as a consequence of Schur's lemma. See Chapter 2B of \cite{DS} for Schur's lemma applications in this context.\\

The first step below follows from Cauchy-Schwarz inequality noting that $\mu_t$ and $U_t$ have support either $A_n$ or its coset depending on $t$. The second follows from Plancherel's formula \eqref{planc} and (\ref{triv}). 
\begin{equation*}
\begin{split}
4 \| \mu_{t} - U_{t} \|_{TV}^2 = \left(\sum\limits_{\sigma \in S_n}|\mu_t(\sigma)-U_t(\sigma)| \right)^2 & \leq \frac{n!}{2} \, \sum\limits_{\sigma \in S_n}|\mu_t(\sigma)-U_t(\sigma)|^2 \\
& = \frac{1}{2} \sum_{\substack{\lambda \, irred. \\ \lambda \neq \lambda_{triv}, \lambda_{sign}}} d_{\lambda} tr(\widehat{\mu_t}(\lambda)\widehat{\mu_t}(\lambda)^{T})
\end{split}
\end{equation*}
\bbox

We evaluate the sum above even further by Lemma \ref{trace} and (\ref{conv}) to arrive at
\begin{equation}\label{ub}
 4\| \mu_{t} - U_{t} \|_{TV}^2 \leq \frac{1}{2}\sum_{\substack{\lambda \, irred. \\ \lambda \neq \lambda_{triv}, \lambda_{sign}}} d_{\lambda}^2 \left(\frac{\chi^{\lambda}_{(n-k, 1^{k})}}{d_{\lambda}}\right)^{2}\left(\frac{\chi^{\lambda}_{(2, 1^{n-2})}}{d_{\lambda}}\right)^{2t}
\end{equation} .

\section{Proof of Theorem \ref{mthm1} } \label{chp2}

\subsection{The moment sequence of $\chi_{\rho}$}

We find the moments of $\chi_{\rho}$ by the formula \eqref{di}. We start with evaluating 
\begin{equation*}
\text{tr}(\widehat{\mu}(\lambda))= d_{\lambda}\left( \frac{\chi^{\lambda}_{(n-1,1)}}{d_{\lambda}} \right)\left( \frac{\chi^{\lambda}_{(2,1^{n-2})}}{d_{\lambda}} \right)^t 
\end{equation*}
for $\lambda \vdash n.$ Observe that for most of the partitions of $n$, the character evaluated at $(n-1)-$cycle gives 0 by the Murnaghan-Nakayama rule (See Chapter 4 of \cite{S} for the Murnaghan-Nakayama rule). Diaconis and Greene \cite{DG} calculated the characters for the remaining partitions, 

\begin{equation*}
\chi^{\lambda}_{(n-1,1)}=
 \left\{
\begin{array}{ll}
      1& \lambda=(n) \\
      (-1)^{|\bar{\bar\lambda}|+1} & \lambda=(n-2-i, 2, 1^{i}) \\
      0 & \text{otherwise} \\
\end{array} 
\right.
\end{equation*}
where $\bar{\bar\lambda}=(\lambda_3,\lambda_4,...,\lambda_m)$.
Furthermore Ding \cite{D} has the following results for the character of $\lambda=(n-2-i, 2, 1^{i})$ evaluated at transpositions,
\begin{equation*}
\frac{\chi^{\lambda}_{(2,1^{n-2})}}{d_{\lambda}}=1-\frac{4+i}{n} 
\end{equation*}
for $i \leq \left\lfloor\frac{n-4}{2}\right\rfloor.$ Therefore we have asymptotically,
\begin{equation} \label{greene}
\text{tr}(\widehat{\mu_{cn+1}}(\lambda))\sim (-1)^{|\bar{\bar\lambda}|+1} e^{-2c(2+|\bar{\bar\lambda}|)}
\end{equation}
for $\lambda=(n-2-i, 2, 1^{i})$ and $i \leq \left\lfloor\frac{n-4}{2}\right\rfloor.$

To compute the moments, we need the following auxillary fact.
\begin{lemma} \label{lemc}
For $x \in \mathbb{R}$, 
\begin{equation*}
\sum_{k=0}^n \frac{1}{k+2} \binom{n}{k} (-x)^{k+2} = \frac{(1-x)^{n+2} -1}{n+2} - \frac{(1-x)^{n+1} -1}{n+1}.
\end{equation*}

\end{lemma}
\pf
Take the derivative of the left hand side with respect to $x.$ Then use binomial theorem to have $x(1-x)^n.$ Next integrate $x(1-x)^n$, and solve for the integral constant considering that the expression on the left hand side is $0$ for $x=0.$ 

\bbox

All needed to calculate the moments of $\chi_{\rho}$ with respect to $\mu_{cn+1}$ are derived. First observe that 
\begin{equation*}
\mathbf{E}_{\mu_{cn+1}}((\chi_{\rho})) = 1,
\end{equation*}
since $\rho=S^{(n)} \oplus S^{(n-1,1)}$ and $\chi^{(n-1,1)}_{(n-1,1)}=0.$ \\

For the higher moments, combine \eqref{coeff} and \eqref{greene} in the formula \ref{di} for $2 \leq r \leq  \left\lfloor\frac{n-4}{2}\right\rfloor$ to arrive at

\begin{align*}
\mathbf{E}_{\mu_{cn+1}}((\chi_{\rho})^r)=& a_{(n),r}+ \sum_{\lambda \vdash n}a_{\lambda,r} tr(\widehat{\mu_{cn+1}}(S^{\lambda})) \\
 \sim & \sum_{i=1}^r {r \brace i} + \sum_{k=|\bar{\bar\lambda}|=0}^{n-4} \sum_{i=|\bar{\lambda}|=k+2}^{r} (k+1) \binom{i}{k+2} {r \brace i} (-1)^{k+1} e^{-2c(k+2)} \\
 =& \sum_{i=1}^r {r \brace i} + \sum_{i=2}^{r} \sum_{k=0}^{i-2} (k+1) \binom{i}{k+2} {r \brace i} (-1)^{k+1} e^{-2c(k+2)} \\
 =& 1+ \sum_{i=2}^{r} {r \brace i} \Big(1-\sum_{i=k+2}^{r} (k+1) \binom{i}{k+2}(-e^{-2c})^{k+2} \Big) \\
 =& 1+ \sum_{i=2}^{r} {r \brace i} \Big(1-i(i-1)\sum_{i=k+2}^{r} \frac{1}{k+2} \binom{i-2}{k}(-e^{-2c})^{k+2} \Big). \\
\end{align*}

Then apply Lemma \ref{lemc} for $x=e^{-2c}$ to have
\begin{equation} \label{momseq}
\begin{split}
\mathbf{E}_{\mu_{cn+1}}((\chi_{\rho})^r)
 \sim & 1+ \sum_{i=2}^{r} {r \brace i} \Big(1-i(i-1) \Big(\frac{(1-e^{-2c})^i-1}{i} - \frac{(1-e^{-2c})^{i-1}-1}{i-1} \Big) \Big) \\
 =& 1 + \sum_{i=2}^{r} {r \brace i} \Big((1-e^{-2c})^i + i e^{-2c} (1-e^{-2c})^{i-1} \Big) \\
 =& \sum_{i=1}^{r} {r \brace i} \Big(1 + i \frac{e^{-2c}}{1-e^{-2c}}\Big)(1-e^{-2c})^{i}. 
\end{split}
\end{equation}

\indent One observation, which is made earlier in \cite{D}, is that the moments obtained above is pretty close to the moments of the Poisson distribution with parameter $(1-e^{-2c}).$ In fact, Kuba and Panholzer \cite{K} shows that if the middle term in the sum above coincided with the $i^{th}$ moment of a distribution, we would have a mixed Poisson random variable where the Poisson paramater has moments proportional to the middle term. But one can easily show that the middle term does not qualify to be a moment sequence. However, Proposition 2 in \cite{K} on the distribution related to the moment sequence of a mixed Poisson distribution can easily be modified to our case, which is in the next section. \\ 

\subsection{The asymptotic distribution of $\chi_{\rho}$}

In this section, we determine the asymptotic distribution of $\chi_{\rho}$ by its moment sequence. We start with a well-known theorem that gives a sufficient condition for the uniqueness of a distribution given its moment sequence. 

\begin{theorem} \emph{(\cite{Bi})} \label{uniq}
Let $X$ be a real random variable having finite moments $\mu_n.$ If the moment generating function $\mathbb{E}(e^{zX})$ of $X$ has positive radius of convergence, then the distribution of $X$ is the only distribution with the moment sequence $\mu_n.$
\end{theorem}

A combinatorial fact that is used for evaluating the sums below is as follows.
\begin{lemma} \emph{(\cite{A})} \label{lem}
If $n>m$,  
\begin{equation*}
\sum_{k=0}^n (-1)^{n-k} \binom{n}{k} k^m = 0.
\end{equation*}
\end{lemma}

Next we verify the hypothesis of Theorem \ref{uniq} for $\chi_{\rho}.$
\begin{proposition}\label{uniq2}
The moment generating function of $\chi_{\rho}$ with distribution $\mu_{cn+1}$ has positive radius of convergence. 
\end{proposition}
\pf 
Let $\zeta \coloneqq 1-e^{-2c}$ and $\eta \coloneqq  \frac{e^{-2c}}{1-e^{-2c}}$ in \eqref{momseq}. Then,
\begin{align*}
\mathbf{E}(e^{z\chi_{\rho}}) =& \sum_{i=0}^{\infty} \mathbf{E}_{\mu_{cn+1}}((\chi_{\rho})^i) \frac{z^i}{i!} \\
=& \sum_{i=0}^{\infty} \sum_{j=0}^{i} {i \brace j} \zeta^j (1 +j \eta) \frac{z^i}{i!} \\
\leq & \sum_{i=0}^{\infty} \sum_{j=0}^{i} {i \brace j} \pi^j (1 +\eta)^j \frac{z^i}{i!} \\
\end{align*}
First applying Lemma \ref{lem}, then using the two-variable generating function identity involving the Stirling numbers of the second kind, which can be found in Chapter 3 of \cite{W}, we obtain

\begin{equation*}
\mathbf{E}(e^{z\chi_{\rho}}) \leq   \sum_{i=0}^{\infty} \sum_{j=0}^{\infty} {i \brace j} \zeta^j (1 +\eta)^j \frac{z^i}{i!} = e^{\zeta (1+ \eta) (e^{z}-1)}.
\end{equation*}

\bbox

Therefore, it is proven by Theorem \ref{uniq} and Proposition \ref{uniq2} that the moment sequence found above uniquely determines the distribution.

The next is a variation of the result discussed above, Proposition 2 in \cite{K}.

\begin{proposition} \label{asdist}
Let $X$ denote a random variable with probability mass function,
\begin{equation*}
P(X=j) = \sum_{i=j}^{\infty} (-1)^{i-j} \binom{i}{j} \alpha_i \frac{\beta^i}{i!}.
\end{equation*}
Then for all $r \in \mathbb{N},$ the moments of $X$ is given by $$\mathbf{E}(X^r)= \sum_{i=1}^{r} {r \brace i} \alpha_i \beta^{i}.$$
\end{proposition}
\pf The $r^{th}$ moment of $X$ can be expressed as
\begin{align*}
\mathbf{E}(X^r)=& \sum_{j=1}^{\infty} \sum_{i=j}^{\infty} (-1)^{i-j} \binom{i}{j} \alpha_i \frac{\beta^i}{i!} j^r \\
=& \sum_{i=1}^{\infty} \sum_{j=1}^{i} (-1)^{i-j} \binom{i}{j} \alpha_i \frac{\beta^i}{i!} j^r. \\
\end{align*}
Then by Lemma \ref{lem},
\begin{align*}
\mathbf{E}(X^r)=& \sum_{i=1}^{\infty} {r \brace i} \alpha_i \frac{\beta^i}{i!}  \\
=& \sum_{i=1}^{r} {r \brace i} \alpha_i \frac{\beta^i}{i!}. 
\end{align*}

\bbox

Finally, we identify the distribution of $\chi_{\rho}$ as $n$ goes to infinity for the symmetric group $S_n.$ If we take $\alpha_i$ to be $\Big(1 + i \frac{e^{-2c}}{1-e^{-2c}}\Big)$ and $\beta$ to be $1-e^{-2c}$ in Proposition \ref{asdist}, we have
\begin{equation} \label{dist}
P({\chi_{\rho}}=j) \sim \sum_{i=j}^{\infty} (-1)^{i-j} \binom{i}{j} \Big(1 + i \frac{e^{-2c}}{1-e^{-2c}}\Big) \frac{(1-e^{-2c})^i}{i!},
\end{equation}
following from Theorem \ref{uniq} and Proposition \ref{uniq2}.
\subsection{Lower bound for the (n-1)-cycle case}

We first calculate the probability that the trace of the defining representation after one $(n-1)-$cycle and $cn$ transpositions is 0. Taking $j=0$ in \eqref{dist},
\begin{align*}
P({\chi_{\rho}}=0) =& \sum_{i=0}^{\infty} (-1)^{i} \binom{i}{0} \Big(1 + i \frac{e^{-2c}}{1-e^{-2c}}\Big) \frac{(1-e^{-2c})^i}{i!} \\
\sim& \sum_{i=0}^{\infty} (-1)^{i} \frac{(1-e^{-2c})^i}{i!} + \sum_{i=0}^{\infty} (-1)^{i}  i\frac{e^{-2c}}{1-e^{-2c}} \frac{(1-e^{-2c})^i}{i!}  \\
=& e^{-(1-e^{-2c})} + \sum_{i=0}^{\infty} (-1)^{i}  i\frac{e^{-2c}}{1-e^{-2c}} \frac{(1-e^{-2c})^i}{i!}  \\
=& e^{-(1-e^{-2c})} - e^{-2c}\sum_{i=0}^{\infty} (-1)^{i} \frac{(1-e^{-2c})^i}{i!} \\
=& e^{-(1-e^{-2c})} - e^{-2c}e^{-(1-e^{-2c})} \\
=& e^{e^{-2c}-1}(1-e^{-2c}).
\end{align*}
Then take $A$ to be the set of fixed point free permutations, since $\chi_{\rho}$ counts the number of fixed points,
\begin{equation}
\mu_{cn+1}(A)=P(\chi_{\rho}=0)\sim \frac{e^{e^{-2c}}(1-e^{-2c})}{e}.
\end{equation} 
Finally plug it in the lower bound inequality \eqref{lb}.
\begin{align*}
\| \mu_{cn+1} - U_{cn+1} \|_{TV} \geq & |\mu_{cn+1}(A) - U_{cn+1}(A)| \\
 \sim & \bigl\lvert \frac{e^{e^{-2c}}(1-e^{-2c})}{e} - \frac{1}{e} \bigr\rvert \\
 = & \frac{ 1-e^{e^{-2c}} + e^{e^{-2c}-2c} }{e}.
\end{align*} 
The proof of Theorem \ref{mthm1} is completed with the result for the upper bound in \cite{D}.

\bbox

\begin{remark}
Unlike the $n$-cycle case, which is studied in \cite{DP}, we have $\mu_{cn+1}(A)< U_{cn+1}(A)$ for the $(n-1)-$cycle case. The comparison is made by the series expansion of the term above, which is
\begin{equation*}
\sum_{i=0}^{\infty} \frac{1}{i+2} \frac{(e^{-2c})^i}{i!} > 0.
\end{equation*}
\end{remark}

\section{Proof of Theorem \ref{mthm2}} \label{chp3}
\subsection{Lower bound for (n-k)-cycle case for $k \geq 2$} \label{lower}

The moment sequence of $\chi_{\rho}$ for any $k \geq 2$ turns out to be quite complicated for identifying the distribution of $\chi_{\rho}$ unlike $k=1$ case above. However, we can derive a lower bound through the first two  moments for $k \geq 3$. The third moment is used only for $k=2$ case. The reason is that the the third moment calculations are cumbersome for $k \geq 3$ and does not yield a significant difference in terms of the lower bound. But for $k=2$, both the calculations are relatively easy and the constant in the lower bound is significantly larger.\\

We have above the decomposition of $\rho$ \eqref{fmom}, and by formula \eqref{rmom} we can find the decomposition of $\rho \otimes \rho$ to be 
\begin{equation}\label{dec}
\rho\otimes \rho =2 S^{(n)} \oplus 3 S^{(n-1,1)} \oplus S^{(n-2,2)} \oplus S^{(n-2,1^2)}.
\end{equation}
Therefore, we have the following expressions for the first two moments.
\begin{equation}\label{mu}
 \begin{split}
\mathbf{E}_{\mu_t}(\chi_{\rho}) &= \text{tr}(\hat{\mu_t}(n))+ \text{tr}(\hat{\mu_t}(n-1,1)) \\
\mathbf{E}_{\mu_t}(\chi_{\rho}^2) &= 2\,\text{tr}(\hat{\mu_t}(n))+ 3 \, \text{tr}(\hat{\mu_t}(n-1,1)) +  \text{tr}(\hat{\mu_t}(n-2,2)) + \text{tr}(\hat{\mu_t}(n-2,1,1))
 \end{split}
\end{equation}
Since $U_t \sim \mathcal{P}(1),$ we also have
\begin{equation}\label{yu}
 \begin{split}
 \mathbf{E}_{U_t}(\chi_{\rho})& \sim 1\\
 \mathbf{E}_{U_t}(\chi_{\rho}^2)& \sim 2\\
 \end{split}
\end{equation}

In order to find the moments needed for the lower bound, we first apply Murnaghan-Nakayama rule to evaluate the characters of the representations in \eqref{dec}. Calculations yield
\begin{equation}\label{karakter}
\begin{split}
\chi^{(n)}_{(n-k,1^k)} &=1 \, \, \text{   for all }k
\\
\chi^{(n-1,1)}_{(n-k,1^k)}&=k-1 \, \, \text{   for all }k\geq 2
\\
\chi^{(n-2,2)}_{(n-k,1^k)}&=
 \left\{
\begin{array}{ll}
      -1 & k=2 \\
       0 & k=3 \\
       1+\frac{(k)(k-3)}{2} & k \geq 4 \\
\end{array} 
\right.
\\
\chi^{(n-2,1^2)}_{(n-k,1^k)}&=
 \left\{
\begin{array}{ll}
       0 & k=2 \\
       1 & k=3 \\
       2+\frac{(k+1)(k-4)}{2} & k \geq 4 \\
\end{array} 
\right.
\end{split}
\end{equation}
The characters on transpositions are evaluated in Chapter 3 of \cite{D}, we present them below and find the asymptotics for $t=\frac{\ln k}{2}n +cn.$
\begin{equation} \label{transpoz}
\begin{split}
\frac{\chi^{(n)}_{(2, 1^{n-2})}}{d_{(n)}} &=1
\\
\left(\frac{\chi^{(n-1,1)}_{(2, 1^{n-2})}}{d_{(n-1,1)}}\right)^t &=  \Big(\frac{n-3}{n-1}\Big)^t \sim \frac{1}{k} e^{-2c}
\\
\left(\frac{\chi^{(n-2,2)}_{(2, 1^{n-2})}}{d_{(n-2,2)}}\right)^t &=  \Big(\frac{n-5}{n-1}\Big)^t \sim \frac{1}{k^2} e^{-4c}
\\
\left(\frac{\chi^{(n-2,1^2)}_{(2, 1^{n-2})}}{d_{(n-2,1^2)}}\right)^t &=  \Big(\frac{n-4}{n}\Big)^t \sim \frac{1}{k^2} e^{-4c}
\end{split}
\end{equation}
So we can write down the moments by evaluating the expressions in \eqref{mu}. 
\begin{equation}\label{tu}
\begin{split}
\mathbf{E}_{\mu_t}(\chi_{\rho}) & \sim 1+ \frac{k-1}{k} e^{-2c} \, \, \text{   for all }k\geq 2 \\
\mathbf{E}_{\mu_t}(\chi_{\rho}^2) & \sim 
 \left\{
\begin{array}{ll}
      2+\frac{3}{2} e^{-2c} - \frac{1}{4} e^{-4c} & k=2 \\
      2+ 2 e^{-2c} + \frac{1}{9} e^{-4c} & k=3 \\
       2+3 \frac{k-1}{k} e^{-2c} + \frac{k^2-3k+1}{k^2} e^{-4c} & k \geq 4 \\
\end{array} 
\right. 
\end{split} 
\end{equation}
Then a lower bound for $k\geq 4$ follows from \eqref{yu} and \eqref{tu} by taking $p=q=2$ in Lemma \ref{p2}, 
\begin{equation}\label{k4}
\begin{split}
 \| \mu_{t} - U_{t} \|_{TV} &\geq \frac{\left(\frac{k-1}{k}e^{-2c}\right)^2}{ 4+3 \frac{k-1}{k} e^{-2c} + \frac{k^2-3k+1}{k^2} e^{-4c}} \\
  &\geq   \frac{\left(\frac{k-1}{k}\right)^2}{4+3 \frac{k-1}{k} + \frac{k^2-3k+1}{k^2}} e^{-4c} \\
  & \geq \frac{3}{35} e^{-4c},  \quad \quad \quad \quad \quad \quad \quad \quad \quad\text{ if  \, $k \geq 4$}
 \end{split}
\end{equation}
We substitute $c=0$ in the denominator for the second inequality, since $c$ can be arbitrarily close to $0$. For the last inequality $k$ is taken $4$, which minimizes the expression. \\
By the same method, we find a lower bound for $k=3.$
 \begin{equation}\label{k3}
\begin{split}
 \| \mu_{t} - U_{t} \|_{TV} &\geq \frac{\left(\frac{2}{3}e^{-2c}\right)^2}{ 4+2 e^{-2c} + \frac{1}{9} e^{-4c}} \\
& \geq \frac{4}{55} e^{-4c},  \quad \quad \quad \quad \quad \quad \quad \quad \quad\text{ if  \, $k =3$}
 \end{split}
\end{equation}
 
For $k=2$, we consider the third moment of $\chi_p$ since the third moment yields a constant twice larger the one obtained by the second moment. We have the decomposition of $\rho^{\otimes 3}$ by \eqref{rmom}, which is
\begin{equation*}
\rho^{\otimes 3}=5 S^{(n)} \oplus 10 S^{(n-1,1)} \oplus 6 S^{(n-2,2)} \oplus 6 S^{(n-2,1^2)}\oplus S^{(n-3,3)}\oplus S^{(n-3,2,1)} \oplus S^{(n-3,1^3)}.
\end{equation*}
By Murnaghan-Nakayama rule, we evaluate the characters of the representations in the decomposition above,
\begin{equation*}
\chi^{\lambda}_{(n-2,1^2)}=
 \left\{
\begin{array}{ll}
      -1 & \lambda=(n-3,3) \\
       0 & \lambda=(n-3,2,1) \\
       0 & \lambda=(n-3,1^3) \\
\end{array} 
\right.
\end{equation*}  
Next we evaluate the non-vanishing character at transpositions and find the asymptotics for $t=\frac{\ln k}{2}n +cn.$
\begin{equation*}
\left(\frac{\chi^{(n-3,3)}_{(2, 1^{n-2})}}{d_{(n-3,3)}}\right)^t =  \Big(\frac{(n-3)(n-4)}{n(n-1)}\Big)^t \sim \frac{1}{k^3} e^{-6c}.
\end{equation*}
Therefore together with \eqref{karakter} and \eqref{transpoz}, we have
\begin{equation*}
\mathbf{E}_{\mu_t}(\chi_{\rho}^3) = 5+ 5 e^{-2c} -\frac{3}{2} e^{-4c} - \frac{1}{8} e^{-6c} 
\end{equation*}
for $k=2.$ Finally apply Lemma \ref{p2} taking $p=3$ and considering that $\mathbf{E}_{U_t}(\chi_{\rho}^2) \sim 5$ to conclude
\begin{equation} \label{k2}
\begin{split}
 \| \mu_{t} - U_{t} \|_{TV} &\geq \frac{\left(\frac{1}{2}e^{-2c}\right)^{3/2}}{ \left(10+ 5 e^{-2c}  -\frac{3}{2} e^{-4c} - \frac{1}{8} e^{-6c} \right)^{1/2}} \\
& \geq \frac{1}{2^{3/2}\sqrt{15}} e^{-3c} \\
& \geq \frac{1}{12} e^{-4c}, 
\quad \quad \quad \quad \quad \quad \quad \quad \quad\text{ if  \, $k =2.$}
 \end{split}
\end{equation}
We compare the bounds in \eqref{k4}, \eqref{k2} and \eqref{k2} to see that $\frac{1}{12} e^{-4c}$ is the smallest, which proves the left hand side of the expression in Theorem \ref{mthm2}.
\begin{remark}
The higher moments can provide better lower bounds. However, the calculations are quite lengthy and it is not essential for our purpose. Therefore the third moment calculations for $k \geq 3$ are skipped.
\end{remark}

\subsection{Auxillary facts for the upper bound}
\begin{figure}[H]
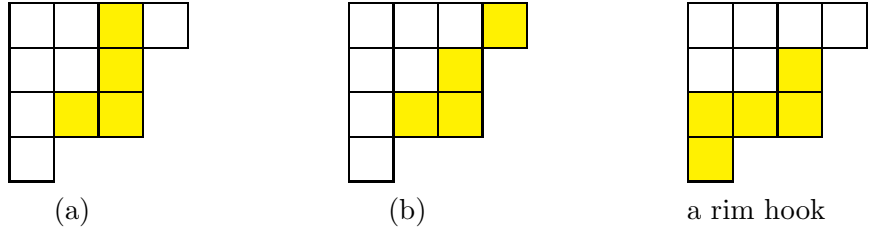
 
\hspace*{1cm}
\ytableaushort{
\empty \empty{*(yellow) \empty}\empty ,\empty \empty{*(yellow) \empty},\empty{*(yellow) \empty}{*(yellow) \empty}, \empty}
\hspace*{2cm}
\ytableaushort{
\empty \empty\empty{*(yellow) \empty} ,\empty \empty {*(yellow) \empty},\empty{*(yellow) \empty}{*(yellow) \empty}, \empty}
\hspace*{2cm}
\ytableaushort{
\empty \empty\empty\empty ,\empty \empty{*(yellow) \empty},{*(yellow) \empty}{*(yellow) \empty}{*(yellow) \empty}, {*(yellow) \empty}}
\newline \newline
\indent \hspace*{1cm} (a) \hspace*{3.7cm} (b) \hspace*{3.2 cm}  a rim hook 
\caption{The shaded region on (a) fails to be a skew diagram. While the shaded region on (b) is a skew diagram, it is not edgewise connected.}\label{rimh}
\end{figure}

The rest of the paper is on the upper bound in Theorem \ref{mthm2}. We first define a central object for the following sections. Then we state two important lemmas which play key role in the evaluation of \eqref{ub}.
\begin{definition}A \textit{skew diagram} of a Young diagram is a subset of cells such that when removed the diagram obtained is a Young diagram. A \textit{rim hook} (also known as \textit{ribbon}) of a Young diagram is an edgewise connected skew diagram with no $2 \times 2$ cells.
\end{definition}
See Figure \ref{rimh} for non-examples and an example of a rim hook.

\begin{lemma} \emph{(\cite{LP})} \label{lem1}
For any Young diagram of shape $\lambda \vdash n$ there exists at most one way to remove a rim hook of length $m > \frac{n}{2}.$
\end{lemma}
A purely combinatorial proof of Lemma \ref{lem1} is found in Section 6 of \cite{LP}.  

\begin{lemma} \emph{(\cite{D})} \label{lem3}
Let $|\lambda|=n$ and $\lambda_1$ be the size of the first row of $\lambda.$ Then,
\begin{equation*}
\sum\limits_{\lambda=(\lambda_1,...)} d_{\lambda}^2 \leq {n \choose \lambda_1}^2 (n- \lambda_1)!.
\end{equation*}
\end{lemma}
The result is stated in Section 3D of \cite{D}. It follows from two results stated earlier in \cite{D}. One is a corollary to Fact 1 in Section 3D, which gives a bound on $d_{\lambda}$. The second one is Corollary 1 in Chapter 2C, which is an important fact in representation theory that the sum of the squares of the dimensions of irreducible representations is equal to the order of the group. 

\subsection{Estimates for the normalized character of transpositions}
In this section, we present estimates for normalized characters, which lead to a monotonicity result that is used to estimate the upper bound on the total variation distance in \ref{UB}. 
\begin{lemma} \emph{(\cite{D})} \label{lem2} Let $\lambda=(\lambda_1, \lambda_2,...,\lambda_m)$ be a partition of $n$ and $(2,1^{n-2})$ stands for transpositions. Then we have the following formula.
\begin{align}
r(\lambda)  \coloneqq \frac{\chi^{\lambda}_{(2, 1^{n-2})}}{d_{\lambda}} =& \frac{1}{n(n-1)} \sum\limits_{i=1}^m \Big(\lambda_i^2 - (2i-1) \lambda_i\Big) \notag \\
 =& \frac{1}{{n \choose 2}} \sum\limits_{i=1}^m {\lambda_i \choose 2}-{\lambda^T_i \choose 2}\label{eqn1}.
\end{align}
\end{lemma}

\begin{lemma}\label{lem6} Let $\lambda=(\lambda_1, \lambda_2,...,\lambda_m)$ be a partition of $n$. Then,
\begin{equation*}
r(\lambda) \leq \frac{\lambda_2-3 + \frac{\lambda_1 (\lambda_1 - \lambda_2 +2)}{n}}{n-1}.
\end{equation*} 
\end{lemma}

\pf By Lemma \ref{lem2},
\begin{equation*}
n(n-1)r(\lambda) =  \sum \lambda_i (\lambda_i -(2i-1))=\lambda_1(\lambda_1-1)+ \lambda_2(\lambda_2-3) + \cdots
\end{equation*}
Then we can bound the expression above as
\begin{align*}
n(n-1) r(\lambda) \leq & \lambda_1 (\lambda_1-1) + (\lambda_2-3) \sum\limits_{i=2} \lambda_i \\
= & \lambda_1 (\lambda_1 - \lambda_2 +2) + \lambda_1(\lambda_2 -3) + (\lambda_2-3) \sum\limits_{i=2} \lambda_i \\
= & (\lambda_2-3)n + \lambda_1 (\lambda_1 - \lambda_2 +2).
\end{align*}

\bbox

We define the following partial order on the set of partitions of $n$, called \textit{dominance order}, which is used to compare the normalized characters of transpositions in the proof of the upper bound. 
\begin{definition} (\cite{M})
If $\lambda=(\lambda_1, \lambda_2,...)$ and $\xi=(\xi_1, \xi_2,...)$ are partitions of $n,$ we say that $\lambda$ \textit{dominates} $\xi$, and denote by $\lambda \succeq \xi,$ if

\begin{equation*}
\sum\limits_{i=1}^j \lambda_i \geq \sum\limits_{i=1}^j \xi_i \textrm{  for all j. }
\end{equation*}
Besides, if $\lambda \succeq \xi$ and $\lambda \neq \xi,$ then we write $\lambda \succ \xi.$
\end{definition}

The following result is in Chapter 3D of \cite{D}.
\begin{lemma} \emph{(\cite{D})}\label{lem5} If $\lambda \succeq \xi,$ then $r(\lambda) \geq r(\xi).$
\end{lemma}
\noindent
An interpretation of the partial order defined above is that $\lambda$ is obtained from $\xi$ by moving boxes up to the right \cite{D}. The following lemma is crucial for evaluating the sums in the next section.

\begin{lemma} \label{lem4}
Suppose $\lambda$ is obtained from $\lambda'$ by moving a box on the rim up to the right. i.e., $$
\lambda'=(\lambda_1, \lambda_2,..., \lambda_{k-1},\lambda_k-1,\lambda_{k+1},..., \lambda_{l-1},\lambda_l+1,\lambda_{l+1},,...) $$
  for some $k<l$ and $\lambda_l < \lambda_k.$  Then $r(\lambda) - r(\lambda')$ is $\frac{2}{n(n-1)}$ times the number of boxes on the rim between the two positions of the displaced box.
\end{lemma}
\pf
\noindent
We calculate the difference in the normalized characters by Lemma \ref{lem2}.
\begin{align*}
r(\lambda) - r(\lambda') =& \frac{1}{n(n-1)}\Bigg( \sum \Big(\lambda_i^2 - (2i-1) \lambda_i\Big) - \Big(\lambda_i^{'2} - (2i-1) \lambda'_i\Big) \Bigg) \\
=& \frac{1}{n(n-1)} \Big( \lambda_k(\lambda_k-(2k-1))-(\lambda_k-1)(\lambda_k-2k)\Big)\\
 + & \frac{1}{n(n-1)} \Big(\lambda_l(\lambda_l-(2l-1) -(\lambda_l+1)(\lambda_l-(2l-2)) \Big)\\
=& \frac{2(\lambda_k-\lambda_l+(l-k)-1)}{n(n-1)}. 
\end{align*}
\bbox

Figure \ref{fig1} below illustrates an example for Lemma \ref{lem4}.

\begin{figure}[H]
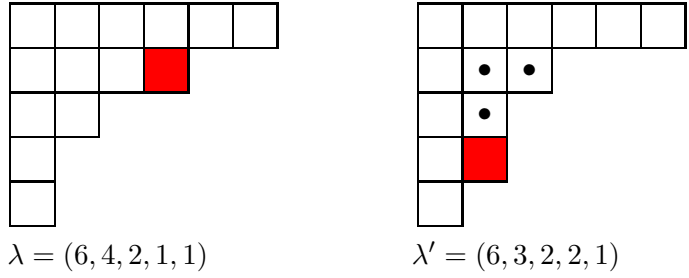

\hspace{1.5cm}
\ytableaushort{
\empty\empty\empty\empty\empty\empty ,\empty\empty\empty{*(red) \empty},\empty\empty,\empty, \empty}
\hspace*{1cm}
\ytableausetup{nobaseline}
\ytableaushort{
\none\empty\empty\empty\empty\empty\empty ,\none\empty\bullet\bullet,\none\empty\bullet,\none\empty{*(red) \empty}, \none\empty} 
\newline \newline
\indent
$  \, \lambda=(6,4,2,1,1) \quad \qquad\qquad\qquad \lambda'=(6,3,2,2,1)$ 
\caption{A pair of Young diagrams that satisfy the hypothesis of Lemma \ref{lem4}.\label{fig1} We have $  r_{\lambda}-r_{\lambda'}=\frac{2}{n(n-1)} \times 3 = \frac{6}{13 \times 12}=\frac{1}{26}.$}
\end{figure}

\begin{remark}
A more general statement for Lemma \ref{lem4}, along with various estimates for normalized characters, is found in \cite{Be}. But Lemma \ref{lem4} is sufficient for our purpose, to limit the growth of successive terms in the upper bound in the next section. 
\end{remark}

\subsection{Upper bound for (n-k)-cycle case for $k \geq 2$} \label{UB}

This final section is devoted to simplifying the sum in \eqref{ub} and bounding it from above. We first express it with the normalized character notation. 
\begin{equation}  \label{sum}
4 \| \mu_{t+1} - U_{t+1} \|^{2}_{TV} \leq  \frac{1}{2} \sum_{\substack{\lambda \vdash n \\ \lambda \neq \lambda_{triv}, \lambda_{sign}}} (\chi^{\lambda}_{(n-k,1^{k})})^2 \, r(\lambda)^{2t}. 
\end{equation}
Consider the first term, $(\chi^{\lambda}_{(n-k,1^{k})})^2$, in the sum. Since $k$ is fixed, eventually $(n-k) > \frac{n}{2}.$ There is thus exactly one way to remove a rim hook of length $(n-k)$ by Lemma \ref{lem1}.  An easy corollary to Murnaghan-Nakayama rule gives 
\begin{equation*}
(\chi^{\lambda}_{(n-k,1^{k})})=\pm d_{\xi}
\end{equation*}
where $\xi$ is the unique diagram obtained from $\lambda$ by removing the rim hook of length $(n-k).$ Then we write the sum \eqref{sum} over the partitions obtained after the removal of a rim hook. \\
 \begin{equation}\label{ubi}
 \begin{gathered}
4 \| \mu_{t+1} - U_{t+1} \|^{2}_{TV} \\ \leq  \frac{1}{2} \left( \sum_{\substack{|\xi|=k \\ \xi \neq (k),(1^{k})}} \sum\limits_{s=1}^{n-k} d_{\xi}^2 \,\, \, r(\lambda^s_{\xi}) ^{2t} +   \sum\limits_{s=2}^{n-k} d_{(k)}^2r(\lambda_{(k)}^s) ^{2t}  +   \sum\limits_{s=1}^{n-k-1} d_{(1^k)}^2r(\lambda_{(1^k)}^s) ^{2t} \right) 
\end{gathered}
\end{equation}
$\lambda_{\xi}^s$ in the sum denotes a partition of $n$ that yields $\xi$ after removing the rim hook, and $s$ is an index for an ordering on the vertical position of the downmost cells of the rim hook among all such partitions. See Figure \ref{fig2} below.  If $n$ is large enough, there are exactly $(n-k)$ possible ways to recover a Young diagram of size $n$ by attaching a rim hook of size $(n-k)$ to a diagram of size $k$. Therefore, $s$ runs through $1$ to $n-k,$ where $\lambda_{\xi}^1=(n-k,1^k).$ The second sum lacks $s=1,$ because it corresponds to the trivial representation. Similarly, the sign representation is excluded in the third sum. Note that if $k=2,$ then  the first tem in the paranthesis above vanishes. \\

\begin{figure}[H]
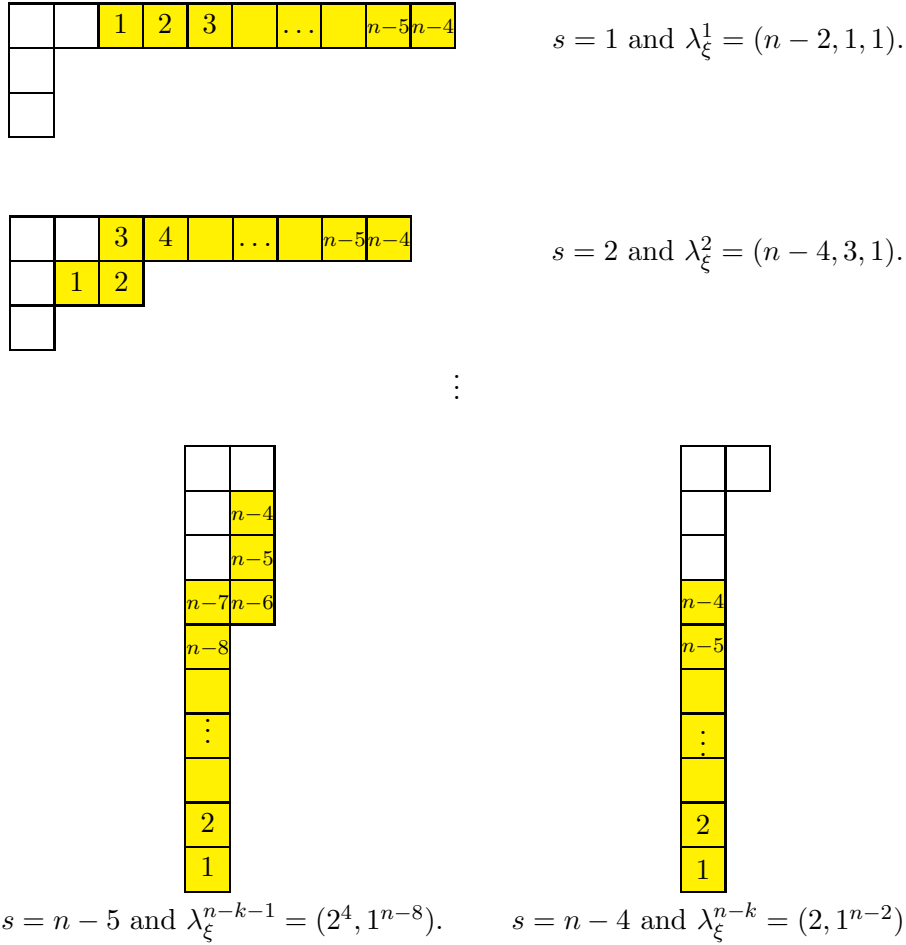

\begin{ytableau}
\empty & \empty & 
 *(yellow){1} & *(yellow) {2} & *(yellow) {3} & *(yellow) \empty   & *(yellow) \dots &*(yellow) \empty & *(yellow) \scriptstyle n-5 & *(yellow) \scriptstyle n-4  \\
\empty \\
\empty \\
\end{ytableau}
\quad\quad\quad $s=1$ and $\lambda_{\xi}^1=(n-2,1,1).$

\vspace*{1cm}

\begin{ytableau}
\empty & \empty &  *(yellow) {3} &
 *(yellow){4} & *(yellow)  \empty   & *(yellow) \dots &*(yellow) \empty & *(yellow)  \scriptstyle n - 5 & *(yellow) \scriptstyle n-4  \\
\empty & 
 *(yellow){1} & *(yellow) {2}    \\
\empty \\
\end{ytableau}
\quad\quad\quad\quad \, $s=2$ and $\lambda_{\xi}^2=(n-4,3,1).$

\begin{center}
$\vdots$
\end{center}

\vspace{5mm} 

\hspace*{1.5cm}
\begin{ytableau}
\empty & \empty   \\
\empty & *(yellow)\scriptstyle n-4  \\
\empty & *(yellow) \scriptstyle n-5\\
*(yellow)\scriptstyle n-7 & *(yellow)\scriptstyle n-6  \\
*(yellow) \scriptstyle n-8 \\  *(yellow)  \empty   \\ *(yellow) \vdots \\
*(yellow) \empty \\ *(yellow) {2} \\
 *(yellow) {1}
\end{ytableau}
\ytableausetup{nobaseline}
\hspace{5cm}
\begin{ytableau}
\empty & \empty  \\
\empty \\
\empty \\
*(yellow)\scriptstyle n-4  \\
*(yellow) \scriptstyle n-5 \\  *(yellow)  \empty   \\ *(yellow) \vdots \\
*(yellow) \empty \\ *(yellow) {2} \\
 *(yellow) {1}
\end{ytableau}
\newline \newline
$s=n-5$ and $\lambda_{\xi}^{n-k-1}=(2^4,1^{n-8}).$ \quad\quad $s=n-4$ and $\lambda_{\xi}^{n-k}=(2,1^{n-2}).$ 
\caption{Young diagrams of $\lambda_{\xi}^1,...,\lambda_{\xi}^{n-k}$ for $k=4$ and $\xi=(2,1,1).$}\label{fig2}
\end{figure}

We regroup the terms in \eqref{ubi} finer. It suffices to consider the partitions with  $\xi_1 \geq \xi^T_1$ where $\xi^T$ is the transpose diagram of $\xi$ and $\xi_1$ is the size of the first row of $\xi.$ The justification comes from a symmetric and an anti-symmetric relation below,
\begin{equation} \label{sym}
\begin{split}
d_{\lambda} &= d_{\lambda^{T}} \\
r(\lambda) &=- r(\lambda^T) \textrm{  by  \eqref{eqn1}  }.
\end{split}
\end{equation}
The second and the third sum in (\ref{ubi}) are the same as well by symmetry  and the fact that $d_{(k)}=d_{(1^k)}=1$. Therefore we have,
\begin{align}
4 \| \mu_{t+1} - U_{t+1} \|^{2}_{TV}  \leq & \frac{1}{2} \, 2 \sum\limits_{i=\ceil[\small]{\sqrt{k}}}^{k-1}  \sum_{\substack{|\xi |=k \notag \\ \xi_1=i}}  \sum\limits_{s=1}^{n-k}  d_{\xi}^2 \, \, \, r(\lambda_{\xi}^s) ^{2t} +  \sum\limits_{s=1}^{n-k-1}  r(\lambda_{\xi}^s) ^{2t}\\
=&   \sum\limits_{i=\ceil[\small]{\sqrt{k}}}^{k-1}  \sum_{\substack{|\xi |=k \\ \xi_1=i}} d_{\xi}^2 \, \sum\limits_{s=1}^{n-k}   r(\lambda_{\xi}^s) ^{2t} +  \sum\limits_{s=1}^{n-k-1}  r(\lambda_{\xi}^s) ^{2t}.\label{eq1}
\end{align}

Next we evaluate the sums $\sum\limits_{s=1}^{n-k}  r(\lambda_{\xi}^s )^{2t}$ and $ \sum\limits_{s=1}^{n-k-1}  r(\lambda_{(k)}^s) ^{2t}$. We identify the largest term in the sum and the difference between the successive terms, and then we bound it from above.
Take $i=\xi_1.$ The first observation is that
\begin{equation*}
\lambda_{\xi}^1 \succ \lambda_{\xi}^2 \succ \cdots \lambda_{\xi}^{n-k-1} \succ \lambda_{\xi}^{n-k}. 
\end{equation*}
So by Lemma \ref{lem5}, the maximum term of the sum is achieved for $s=1$. The bound for that term is derived from Lemma \ref{lem6} below. 
\begin{equation} \label{obs1}
\begin{split}
r(\lambda_{\xi}^1)  \leq &  \frac{i-3 + \frac{(n-k +i) (n - k +2)}{n}}{n-1} \\
=& \frac{n-2k+2i-1 + o(1)}{n-1} \\
=& 1 - \frac{2(k-i)+o(1)}{n-1}
\end{split}
\end{equation}
 The second observation is that the last term in the sum in absolute value admits the same upper bound with $r(\lambda_{\xi}^1)$, i.e., 
\begin{equation} \label{obs2}
\begin{split}
-r(\lambda_{\xi}^{n-k})=r((\lambda_{\xi}^{n-k})^T) \leq & \frac{i-3 + \frac{(n-k+i) (n-k+2)}{n}}{n-1} \\  = & 1 - \frac{2(k-i)+o(1)}{n-1}.
\end{split}
\end{equation}
by \eqref{sym} \\ \\
\noindent Finally, it follows from Lemma \ref{lem4} that for the successive diagrams we have
\begin{equation}\label{obs3}
r(\lambda_{\xi}^i)-r(\lambda_{\xi}^{i+1})\geq \frac{2(n-k-1)}{n(n-1)} = \frac{2(n-k)-2}{n(n-1)} = \frac{2-o(1)}{n-1}.
\end{equation}
Lemma \ref{lem4} is applicable above. Because $\lambda^{i+1}$ can be thought to be obtained from $\lambda^i$ by moving the box on the rightmost and uppermost end of the rim hook of length $(n-k)$ to the leftmost and downmost end, if it results in another Young diagram. If not, it is moved even further left or further below. Therefore the number of boxes on the rim between two positions is greater than or equal to $n-k-1$. \\ \\
\noindent We have thus identified the largest term in the sum and the difference between the successive terms.To bound the sum, first take $p$ to be the largest index for which $r(\lambda_p)$ is positive in order to split the inner sum for negative and positive normalized characters; the justification is by the symmetry \eqref{obs2} observed above. Then by \eqref{obs1} and \eqref{obs3}, we have
\begin{equation}\label{normzero}
\sum\limits_{s=1}^p   r(\lambda_{\xi}^s )^{2t} \leq  \sum\limits_{s=0}^{p-1} \Big( 1 - \frac{2(k-i)+2s+o(1)}{n-1}\Big)^{2t} 
\end{equation}
Taking $t=cn + \frac{n}{2} \ln k$ in \eqref{normzero},
\begin{equation}\label{norm}
\begin{split}
\sum\limits_{s=1}^p   r(\lambda_{\xi}^s )^{2t} &\leq  \sum\limits_{s=0}^{p-1} \Big( 1 - \frac{2(k-i)+2s+o(1)}{n-1}\Big)^{2cn + n \ln k}  \\
& \sim  \sum\limits_{s=0}^{p-1} e^{-4(k-i+s)c} k^{-2(k-i+s)}. 
\end{split}
\end{equation}

Similarly,
\begin{equation}\label{norm2}
\begin{split}
\sum\limits_{s=p+1}^{n-k}   r(\lambda_{\xi}^s )^{2t} &\leq  \sum\limits_{s=0}^{n-k-p-1} \Big( 1 - \frac{2(k-i)+2s+o(1)}{n-1}\Big)^{2t} \\
& \sim  \sum\limits_{s=0}^{p-1} e^{-4(k-i+s)c} k^{-2(k-i+s)}.
\end{split}
\end{equation}

We add \eqref{norm} and \eqref{norm2} below to arrive at an upper bound for the sum. 
\begin{align*}
\sum\limits_{s=1}^{n-k}  r(\lambda_{\xi}^s )^{2t} \leq & 2 \sum\limits_{s=0}^{\infty} e^{-4(k-i+s)c} k^{-2(k-i+s)} \\
 =& 2e^{-4(k-i)c} k^{-2(k-i)}\sum\limits_{s=0}^{\infty} (e^{-4c}k^{-2})^s \\
 =&  \frac{2}{1-e^{-4c}k^{-2}}e^{-4(k-i)c} k^{-2(k-i)}  .
\end{align*}

The argument above also works for  $\sum\limits_{s=2}^{n-k} r(\lambda_{(k)}^s) ^{2t}$ , the second term in \eqref{eqn1}. Consider
\begin{equation*}
\lambda_{(k)}^2 \succ \lambda_{(k)}^3 \succ \cdots \lambda_{(k)}^{n-k-1} \succ \lambda_{(k)}^{n-k}, 
\end{equation*}
where $\lambda^2_{(k)}=(n-k-1,k+1)$ and $\lambda^{n-k}_{(k)}=(k,1^{n-k})$. Apply Lemma \ref{lem6} to have
\begin{equation*}
\begin{split}
r(\lambda_{(k)}^2) & \leq  1- \frac{2k+2+ o(1)}{n-1}, \\
-r(\lambda_{(k)}^{n-2}) & \leq  1- \frac{2k-2+ o(1)}{n-1}. 
\end{split}
\end{equation*}
So the leading term in absolute value is bounded by $1 - \frac{2(k-i)+o(1)}{n-1}$. We repeat \eqref{norm} and \eqref{norm2}, and conclude
\begin{equation*}
\sum\limits_{s=2}^{n-k}  r(\lambda_{(k)}^s) ^{2t} \leq \frac{2}{1-e^{-4c}k^{-2}}e^{-4(k-1)c} k^{-2(k-1)}.
\end{equation*}

Going back to the upper bound inequality \eqref{eq1}, we replace the sum of normalized characters in the expression by the bounds found for them above.
\begin{align*}
 4 & \| \mu_{t+1} - U_{t+1} \|^{2}_{TV} \\ \leq & \frac{2}{1-e^{-4c}k^{-2}} \left(\sum\limits_{i=\ceil[\small]{\sqrt{k}}}^{k-1}  e^{-4(k-i)c} k^{-2(k-i)}\sum_{\substack{|\xi |=k \\ \xi_1=i}} d_{\xi}^2 \, +  e^{-4(k-1)c} k^{-2(k-1)} \right) \\
 \leq & \frac{2}{1-e^{-4c}k^{-2}} \left(\sum\limits_{i=\ceil[\small]{\sqrt{k}}}^{k-1} {k \choose i}^2 (k - i)!  e^{-4(k-i)c} k^{-2(k-i)} +  e^{-4(k-1)c} k^{-2(k-1)} \right) \\
\end{align*}
The second inequality follows from Lemma \ref{lem3}. We change the index from $i$ to $j= k-i$ for simplicity, and simplify the expression even further. For $k \geq 2$ and $c > 0$,  
\begin{align*}
 4 & \| \mu_{t+1} - U_{t+1} \|^{2}_{TV} \\ \leq & \frac{2}{1-e^{-4c}k^{-2}} \left( \sum\limits_{j=1}^{\floor[\small]{k-\sqrt{k}}}  {k \choose j}^2 j! \, e^{-4jc} k^{-2j} +  e^{-4(k-1)c} k^{-2(k-1)} \right)                    
\\
 = &  \frac{2}{1-e^{-4c}k^{-2}}   \left( \sum\limits_{j=1}^{\floor[\small]{k-\sqrt{k}}}  \frac{1}{j!} \, \frac{k^2  \cdots (k-j+1)^2 }{k^{2j}} \, e^{-4(j-1)c} + e^{-4(k-2)c} k^{-2(k-1)}\right)e^{-4c} \\
  \leq &  \frac{2}{1-k^{-2}}   \left( \sum\limits_{j=1}^{\floor[\small]{k-\sqrt{k}}}  \frac{1}{j!} \, \frac{k^2  \cdots (k-j+1)^2 }{k^{2j}} \,  +  k^{-2(k-1)}\right)e^{-4c}\\
   \leq & 2 \left( \sum\limits_{j=1}^{\infty} \frac{1}{j!} \right) e^{-4c}\\
     = & 2 (e-1) e^{-4c}.
\end{align*}
This, combined with the lower bound result in \ref{lower}, completes the proof of Theorem \ref{mthm2}.

\bbox


\begin{thebibliography}{1}

\bibitem{A} Aigner, M., A Course in Enumeration, Springer, Berlin, Germany (2007)

\bibitem{B} Behrends, E., Introduction to Markov Chains (with Special Emphasis on
Rapid Mixing), Vieweg Verlag, Braunschweig/Wiesbaden (2000)

\bibitem{BSZ} Berestycki, N., Schramm, O., Zeitouni, O., Mixing times for random k-cycles and coalescence-fragmentation chains, The Annals of Probability \textbf{39}, No. 5, 1815-1843 (2011)

\bibitem{Be} Bernstein, M., Likelihood Orders for the $p$-Cycle Walks on the Symmetric Group, The Electronic Journal of Combinatorics \textbf{25}, No. 1, 1-25 (2018)

\bibitem{Bi} Billingsley, P., Probability and Measure,  Second Edition , p.406. John Wiley and Sons, New York, NY (1986)

\bibitem{Bo} Bormashenko, O., A coupling argument for the random transposition walk, arXiv preprint arXiv:1109.3915 (2011)

\bibitem{D} Diaconis, P., Group Representations in Probability and Statistics, Institute of Mathematical Sciences, Lecture Notes-Monograph Series \textbf{11}, Hayward, CA (1988)

\bibitem{DG} Diaconis, P., and Greene, C., Applications of Murphy's elements, Stanford University Technical Reports No.335, 1-22 (1989)  

\bibitem{DS} Diaconis, P., Shahshahani, M., Generating a random permutation with random transpositions. Z, Wahrsch. Verw. Gebiete \textbf{57} No. 2 159-179 (1981)

\bibitem{DP} Ding, S., A Random Walk in Representations, Ph.D. dissertation, University of Pennslyvania (2014)


\bibitem{GC} Goupil, A., Chauve C., Combinatorial operators for Kronecker powers
of representations of $S_n$, S\'eminaire Lotharingiene de Combinatoire, \textbf{54}, B54j (2006)

\bibitem{J} James G. D., The Representation Theory of the Symmetric Groups, p.8. Springer-Verlag, Berlin (1978)

\bibitem{K} Kuba, M., Panholzer, A., On Moment Sequences and Mixed Poisson Distribution, Probability Surveys \textbf{13}, 89-155 (2016) 

\bibitem{L} Levin, D.A., Peres, Y., Wilmer, E.L., Markov Chains and Mixing Time, AMS, Providence, Rhode Island (2009)

\bibitem{LP} Lulov, N., Pak, I., Rapidly Mixing Random Walks and Bounds on Characters of the Symmetric Group, Journal of Algebraic Combinatorics \textbf{16} No. 2 , 151-163(2002)

\bibitem{M} Macdonald, I. G., Symmetric functions and Hall polynomials, Oxford University Press, New York, NY (1979)

\bibitem{Sa} Sagan, B., The Symmetric Group, Brooks/Cole Publishing Co., Belmont, California, 1991. 

\bibitem{S} Saloff-Coste, L., Random walks on finite groups, Probability on Discrete Structures, Encyclopaedia
Math. Sci.\textbf{110}, 263-346 (2004) 

\bibitem{T}Tak\'{a}cs, L., The problem of coincidences, Archive for History of Exact Sciences \textbf{21}, No. 3, 229-244 (1980)

\bibitem{W} Wilf, H. S., Generatingfunctionology, Academic Press, New York, NY (1990)  

\end{thebibliography}
\end{document}